\documentclass{article}
\usepackage{t1enc}
\usepackage[latin1]{inputenc}
\usepackage[english]{babel}
\usepackage{amsmath,amsthm}
\usepackage{amsfonts}
\usepackage{latexsym}
\usepackage[dvips]{graphicx}
\usepackage{graphicx}
\newtheorem{theorem}{Theorem}[section]

\newtheorem{corollary}[theorem]{Corollary}

\begin{document}

\title{\textbf{On Potentially $(K_5-H)$-graphic
Sequences}}

\author{Lili Hu  and  Chunhui Lai \thanks{Rearch is supported by NNSF of China(10271105)
and by NSF of Fujian(Z0511034), Fujian Provincial Training Foundation for
"Bai-Quan-Wan Talents Engineering", Project of Fujian Education
Department and Project of Zhangzhou Teachers College.}
\\{\footnotesize Department of Mathematics}
\\{\footnotesize Zhangzhou Teachers College, Zhangzhou}
\\{\footnotesize  Fujian 363000, P. R. of China}
\\{\footnotesize  jackey2591924@163.com(Lili Hu)}
\\{\footnotesize  zjlaichu@public.zzptt.fj.cn(Chunhui
 Lai, Corresponding author)}
\and Ping Wang \thanks{Research is support by NSERC}
\\{\footnotesize Dept. of Math., Stats. and Computer Science}\\
{\footnotesize St. Francis Xavier University}\\
{\footnotesize Antigonish, NS, Canada, B2G 2W5}
\\{\footnotesize pwang@stfx.ca(Ping Wang)}}

\date{}

\maketitle
\begin{abstract}
Let $K_m-H$ be the graph obtained from $K_m$ by removing the edges
set $E(H)$ of $H$ where $H$ is a subgraph of $K_m$. In this paper,
we characterize the potentially  $K_5-P_4$ and $K_5-Y_4$-graphic
sequences where $Y_4$ is a tree on 5 vertices and 3 leaves.
\end{abstract}

\section{Introduction}
We consider finite simple graphs. Any undefined notation follows
that of Bondy and Murty $[1]$. An $n$-term non-increasing
nonnegative integer sequence $\pi=(d_1,d_2,\cdots,d_n)$ is said to
be graphic if it is the degree sequence of a simple graph $G$ of
order $n$; such a graph $G$ is referred as a realization of $\pi$.
Let $\sigma(\pi)$ the sum of all the terms of $\pi$, and let $[x]$
be the largest integer less than or equal to $x$.  Let $Y_4$
denote a tree on 5 vertices and 3 leaves. A graphic sequence $\pi$
is said to be potentially $H$-graphic if it has a realization $G$
containing $H$ as a subgraph. Let $G-H$ denote the graph obtained
from $G$ by removing the edges set $E(H)$ where $H$ is a subgraph
of $G$.  In the degree sequence, $r^t$ means $r$ repeats $t$
times, that is, in the realization of the sequence there are $t$
vertices of degree $r$.

In 1907, Mantel first proposed the problem of determining the
maximum number of edges in a graph without containing $3$-cycles. In
general, this problem can be phased as determining the maximum
number of edges, denoted $ex(n,H)$, of a graph with $n$ vertices not
containing $H$ as a subgraph. This area of research is called
extremal graph theory. In terms of graphic sequences, the number
$2ex(n,H)+2$ is the minimum even integer $l$ such that every
$n$-term graphical sequence $\pi$ with $\sigma (\pi)\geq l $ is
forcibly $H$-graphical. In 1991, Erd\"os,\ Jacobson and Lehel $[2]$
showed $\sigma(K_k, n) \ge (k-2)(2n-k+1)+2$ and conjectured that the
equality holds. In the same paper, they proved the conjecture is
true for the case $k=3$ and $n \ge 6$. The cases $k=4$ and $5$ were
proved separately in $[3]$, $[11]$ and $[12]$. Based on linear
algebraic techniques, Li, Song and Luo $[13]$ proved the conjecture
true for $k \geq 6$ and $n \geq \binom{k}{2}+3$. Recently, Ferrara,
Gould and Schmitt  proved the conjecture $[4]$ and they also
determined in $[5]$ $\sigma(F_k,n)$ where $F_k$ denotes the graph of
$k$ triangles intersecting at exactly one common vertex.

In 1999, Gould, Jacobson and Lehel $[3]$ considered the following
generalized problem: determine the smallest even integer
$\sigma(H,n)$ such that every $n$-term positive graphic sequence
$\pi=(d_1,d_2,\cdots,d_n)$ with $\sigma(\pi)\geq \sigma(H,n)$ has a
realization $G$ containing $H$ as a subgraph. They proved
$\sigma(pK_2, n)=(p-1)(2n-p)+2$ for $p\ge 2$ and $\sigma(C_4,
n)=2[\frac{3n-1}{2}]$ for $n\ge 4$. Lai $[8]$ determined
$\sigma(K_4-e,n)$ for $n\geq4$ . Yin, Li, and Mao $[19]$ determined
$\sigma(K_{r+1}-e,n)$ for $r\geq3$ and $r+1\leq n\leq 2r$ and
$\sigma(K_5-e,n)$ for $n\geq5$, and Yin and Li $[18]$ further
determined $\sigma(K_{r+1}-e,n)$ for $r\geq2$ and $n\geq3r^2-r-1$.
Moreover, Yin and Li in $[18]$ also gave two sufficient conditions
for a sequence $\pi\epsilon GS_n$ to be potentially
$(K_{r+1}-e)$-graphic. Yin $[21]$ determined $\sigma(K_{r+1}-K_3,n)$
for $r\geq3$ and $n\geq 3r+5$.  Lai $[9]$ determined
$\sigma(K_5-P_3,n)$ and $\sigma(K_5-P_4,n)$ for $n\geq5$. Lai and Hu
$[10]$ determined $\sigma(K_{r+1}-H,n)$ for $n\geq4r+10$, $r\geq3$,
$r+1\geq k\geq4$ and $H$ be a graph on $k$ vertices which containing
a tree on 4 vertices but not containing a cycle on 3 vertices and
$\sigma(K_{r+1}-P_2,n)$ for $n\geq4r+8$, $r\geq3$.

A harder question is to characterize the potentially $H$-graphic
sequences without zero terms.  That is, finding necessary and
sufficient conditions for a sequence to be a $H$-graphic sequence.
Luo $[15]$ characterized the potentially $C_k$-graphic sequences
for each $k=3,4$ and $5$. Recently, in $[16]$, Luo and Warner also
characterized the potentially $K_4$-graphic sequences. Eschen and
Niu $[17]$ characterized the potentially $(K_4-e)$-graphic
sequences. Hu and Lai $[6]$ characterized the potentially
$(K_5-C_4)$-graphic sequences. Yin and Chen $[20]$ characterized
the
 potentially $K_{r,s}$-graphic sequences for $r=2,s=3$ and
 $r=2,s=4$, where $K_{r,s}$ is an $r\times s$ complete bipartite
 graph.

In attempt to completely characterize the potentially $K_5-H$ -
graphic sequences, we will characterize the potentially $K_5-P_4$
and $K_5-Y_4$ - graphic sequences in this paper.

Let $\pi=(d_1,d_2,\cdots,d_n)$ be a nonincreasing positive integer
sequence. We write $m(\pi)$ and $h(\pi)$ to denote the largest
positive terms of $\pi$ and the smallest positive terms of $\pi$,
respectively. $\pi^{\prime\prime} =(d_1-1,d_2-1, \cdots,
d_{d_{n}}-1,d_{{d_n}+1}, \cdots, d_{n-1})$ is the residual sequence
obtained by laying off $d_n$ from $\pi$. We denote
$\pi^\prime=(d_1^\prime, d_2^\prime, \cdots, d_{n-1}^\prime)$ where
$d_1^\prime\geq d_2^\prime\geq \cdots \geq d_{n-1}^\prime$ is a
rearrangement of the $n-1$ terms in $\pi^{\prime\prime}$. We denote
$\pi^\prime$ the residual sequence obtained by laying off $d_n$ from
$\pi$ and all the graphic sequences have no zero terms. We need the
following results.

\begin{theorem}$[3]$\label{could}
If $\pi=(d_1,d_2,\cdots,d_n)$ is a graphic sequence with a
realization $G$ containing $H$ as a subgraph, then there exists a
realization $G^\prime$ of $\pi$ containing $H$ as a subgraph so
that the vertices of $H$ have the largest degrees of $\pi$.
\end{theorem}

\begin{theorem}$[14]$
If $\pi=(d_1,d_2,\cdots,d_n)$ is a sequence of nonnegative integers
with $1\leq m(\pi)\leq2$, $h(\pi)=1$ and even $\sigma(\pi)$, then
$\pi$ is graphic.
\end{theorem}

\begin{theorem}$[7]$
$\pi$ is graphic if and only if $\pi^\prime$ is graphic.
\end{theorem}

The following corollary is obvious.

\begin{corollary}\label{coro}
Let $H$ be a simple graph. If $\pi^\prime$ is potentially
$H$-graphic, then $\pi$ is potentially $H$-graphic.
\end{corollary}

\section{Main Theorems}

\begin{theorem}\label{the1}
Let $\pi=(d_1,d_2,\cdots,d_n)$ be a graphic sequence with $n\geq5$.
Then $\pi$ is potentially $(K_5-P_4)$-graphic if and only if the
following conditions hold:
\begin{enumerate}
\item $d_2\geq3$.

\item $d_5\geq2$.

\item $\pi\neq(n-1,k,2^t,1^{n-2-t})$ where $n\geq5$,
$k,t=3,4,\cdots,n-2$, and, $k$ and $t$ have different parities.

\item For $n\geq5$, $\pi\neq(n-k,k+i,2^i,1^{n-i-2})$ where
$i=3,4,\cdots,n-2k$ and $k=1,2,\cdots,[\frac{n-1}{2}]-1$.

\item If $n=6,7$, then $\pi\neq (3^2,2^{n-2})$.
\end{enumerate}

\end{theorem}

\begin{proof}

First we show the conditions (1)-(5) are necessary conditions for
$\pi$ to be potentially $(K_5-P_4)$-graphic. Assume that $\pi$ is
potentially $(K_5-P_4)$-graphic. $(1)$, $(2)$ and $(5)$ are obvious.
If $\pi=(n-1,k,2^t,1^{n-2-t})$ is potentially $(K_5-P_4)$-graphic,
then according to Theorem 1.1, there exists a realization $G$ of
$\pi$ containing $K_5-P_4$ as a subgraph so that the vertices of
$K_5-P_4$ have the largest degrees of $\pi$. Therefore, the sequence
$\pi^*=(n-4,k-3,2^{t-3},1^{n-2-t})$ obtained from $G-(K_5-P_4)$ must
be graphic. Since the edge between two vertices with degree $n-4$
and $k-3$ has been removed from the realization of $\pi^*$, thus,
$\Delta(G-(K_5-P_4))\leq n-5$, a contradiction. Hence, $(3)$ holds.
If $\pi=(n-k,k+i,2^i,1^{n-i-2})$ is potentially $(K_5-P_4)$-graphic,
then according to Theorem 1.1, there exists a realization $G$ of
$\pi$ containing $K_5-P_4$ as a subgraph so that the vertices of
$K_5-P_4$ have the largest degrees of $\pi$. Therefore, the sequence
$\pi^*=(n-k-3,k+i-3,2^{i-3},1^{n-i-2})$ obtained from $G-(K_5-P_4)$
must be graphic and there is no edge between two vertices with
degree $n-k-3$ and $k+i-3$ in the realization of $\pi^*$. Let $G^*$
be a realization of $\pi^*$, and, $d_{G^*}(x)=n-k-3$ and
$d_{G^*}(y)=k+i-3$. Consider a partition of $G^*$ where $X=\{x,y\}$
and $Y=V(G^*)-\{x,y\}$. It follows that the number of edges between
$X$ and $Y$ equals $(n-k-3)+(k+i-3)\leq2(i-3)+(n-i-2)$, that is,
$[(n-k-3)+(k+i-3)]-[2(i-3)+(n-i-2)]=2 \leq 0$, a contradiction.
Hence, $(4)$ holds.

Now we turn to show the conditions (1)-(5) are sufficient
conditions for $\pi$ to be potentially $(K_5-P_4)$-graphic.
Suppose the graphic sequence $\pi$ satisfies the conditions (1) to
(5). Our proof is by induction on $n$. We first prove the base
case where $n=5$. Since $\pi\neq(4^2,2^3)$, then $\pi$ must be one
of the following sequences: $(4^5)$, $(4^3,3^2)$, $(4^2,3^2,2)$,
$(4,3^4)$, $(4,3^2,2^2)$, $(3^4,2)$, $(3^2,2^3)$. It is easy to
check that all of these are potentially $(K_5-P_4)$-graphic. Now
we assume that the sufficiency holds for $n-1$ ($n\geq6$). We will
prove $\pi$ is
potentially $(K_5-P_4)$-graphic.\\

\noindent \textbf{Case 1:  $\pi^\prime=(3^2,2^4)$}

Clearly, $n=7$ and $\pi$ must be one the following sequences
$(4^2,2^5)$, $(4,3^2,2^4)$, $(3^4,2^3)$, $(4,3,2^4,1)$ or
$(3^3,2^3,1)$. It is easy to check that
all of these are potentially $(K_5-P_4)$-graphic.\\

\noindent\textbf{Case 2:  $\pi^\prime=(3^2,2^5)$}

Clearly, $n=8$ and $\pi$ must be one the following sequences
$(4^2,2^6)$, $(4,3^2,2^5)$, $(3^4,2^4)$, $(4,3,2^5,1)$ or
$(3^3,2^4,1)$. It is easy to check that all of these are
potentially
$(K_5-P_4)$-graphic.\\

\noindent\textbf{Case 3: $d_n\geq3$}

Clearly, $\pi^\prime$ satisfies the assumption, and thus, by the
induction hypothesis, $\pi^\prime$ is potentially
$(K_5-P_4)$-graphic, and hence so is $\pi$. In the following, we
only consider the cases $d_n=1$ or $d_n=2$.\\

\noindent\textbf{Case 4:} $\pi^\prime=(n-2,k,2^t,1^{n-3-t})$ where
$n-1\geq5$, $k,t=3,4,\cdots,n-3$, and, $k$ and $t$ have different
parities

If $d_n=2$, then $\pi^{\prime}=(n-2,k,2^{n-3})$. If $k\geq4$, then
$\pi=(n-1,k+1,2^{n-2})$ which contradicts condition (3). If $k=3$,
that is $\pi^{\prime}=(n-2,3,2^{n-3})$, then $\pi=(n-1,4,2^{n-2})$
or $\pi=(n-1,3^2,2^{n-3})$. But $\pi=(n-1,4,2^{n-2})$ contradicts
condition (3), thus $\pi=(n-1,3^2,2^{n-3})$ where $n$ is odd. We
will show $\pi=(n-1,3^2,2^{n-3})$ is potentially
$(K_5-P_4)$-graphic. In other words, we would like to show
$\pi_1=(n-4,2^{n-5},1)$ is graphic. It suffices to show
$\pi_2=(1^{n-5})$ where $n \geq 7$ is graphic. By $\sigma(\pi_2)$
being even and Theorem 1.2, $\pi_2$ is graphic. Thus,
$\pi=(n-1,3^2,2^{n-3})$ is potentially $(K_5-P_4)$-graphic.

If $d_n=1$, then $\pi=(n-1,k,2^t,1^{n-2-t})$ which contradicts
condition (3).\\

\noindent\textbf{Case 5:} $\pi^{\prime}=(n-1-k,k+i,2^i,1^{n-i-3})$
where $i=3,4,\cdots,n-1-2k$ and $k=1,2,\cdots,[\frac{n}{2}]-2$.
\par
 If $d_n=2$,
then $n-i-3=0$ and $\pi=(n-k,k+i+1,2^{i+1})$ and it contradicts
condition (4). \par
  If $d_n=1$ and $n-1-k=k+i+1$, then $\pi=(n-k,k+i,2^i,1^{n-i-2})$ or $\pi=((n-1-k)^2,2^i,1^{n-i-2})$, both of these are
 contradict to condition (4). If $d_n=1$ and $n-1-k=k+i$ or $n-1-k\geq k+i+2$, then
$\pi=(n-k,k+i,2^i,1^{n-i-2})$ which also contradicts condition (4).\\

\noindent\textbf{Case 6:} $d_n=2$,
$\pi^\prime\neq(n-2,k,2^{n-3})$,
$\pi^\prime\neq(n-1-k,n+k-3,2^{n-3})$, $\pi^\prime\neq(3^2,2^4)$,
and $\pi^\prime\neq(3^2,2^5)$

Consider
$\pi^\prime=(d_1^\prime,d_2^\prime,\cdots,d_{n-1}^\prime)$. Since
$d_2\geq3$, we have $d_{n-1}^\prime\geq2$. If $d_2^\prime\geq 3$,
then $\pi^\prime$ satisfies the assumption. Thus, $\pi^\prime$ is
potentially $(K_5-P_4)$-graphic. Hence, we may assume
$d_2^\prime=2$, that is, $d_2=3$ and $d_3=d_4=\cdots=d_n=2$. It
follows $\pi=(d_1,3,2^{n-2})$. Since $\sigma(\pi)$ is even, $d_1$
must be odd. If $d_1=3$, then $\pi=(3^2,2^{n-2})$. Since
$\pi\neq(3^2,2^4)$ and $\pi\neq(3^2,2^5)$, we have $n\geq8$. We
will show $\pi$ is potentially $(K_5-P_4)$-graphic. It suffices to
show $\pi_1=(2^{n-5})$ is graphic. Clearly, $C_{n-5}$ is a
realization of $\pi_1$. If $d_1\geq5$, since
$\pi\neq(n-1,3,2^{n-2})$, we have $d_1\leq n-2$. We will prove
$\pi=(d_1,3,2^{n-2})$ where $d_1\geq5$ and $n\geq d_1+2$ is
potentially $(K_5-P_4)$-graphic. We would like to show
$\pi_1=(d_1-3,2^{n-5})$ is graphic. It suffices to show
$\pi_2=(2^{n-d_1-2},1^{d_1-3})$ is graphic. Since $\sigma(\pi_2)$
is even, $\pi_2$ is graphic by Theorem 1.2.
Thus, $\pi=(d_1,3,2^{n-2})$ is potentially $(K_5-P_4)$-graphic.\\

\noindent\textbf{Case 7:} $d_n=1$,
$\pi^\prime\neq(n-2,k,2^t,1^{n-3-t})$,
$\pi^\prime\neq(n-1-k,k+i,2^i,1^{n-i-3})$,
$\pi^\prime\neq(3^2,2^4)$, and $\pi^\prime\neq(3^2,2^5)$

Consider
$\pi^\prime=(d_1^\prime,d_2^\prime,\cdots,d_{n-1}^\prime)$. Since
$d_2\geq3$ and $d_5\geq2$, we have $d_1^\prime\geq3$ and
$d_5^\prime\geq2$. If $d_2^\prime\geq3$, then $\pi^\prime$
satisfies the assumption. Thus, $\pi^\prime$ is potentially
$(K_5-P_4)$-graphic. Hence, we may assume $d_2^\prime=2$, that is,
$d_1=d_2=3$ and $d_3=d_4=d_5=2$. Thus $\pi=(3^2,2^t,1^{n-2-t})$
where $t\geq3$ and $n-2-t\geq1$. Since $\sigma(\pi)$ is even,
$n-2-t$ must be even. We will prove $\pi$ is potentially
$(K_5-P_4)$-graphic. It suffices to show
$\pi_1=(2^{t-3},1^{n-2-t})$ is graphic.  Since $\sigma(\pi_1)$ is
even, $\pi_1$ is graphic by Theorem 1.2 and, in turn, $\pi$ is
potentially $(K_5-P_4)$-graphic.

This completes the proof.
\end{proof}

\begin{theorem}\label{the2}
Let $\pi=(d_1,d_2,\cdots,d_n)$ be a graphic sequence with $n\geq5$.
Then $\pi$ is potentially $(K_5-Y_4)$-graphic if and only if the
following conditions hold:
\begin{enumerate}
\item $d_3\geq3$.

\item $d_4\geq2$.

\item $\pi\neq(3^6)$.
\end{enumerate}
\end{theorem}

\begin{proof}
Assume that $\pi$ is potentially $(K_5-Y_4)$-graphic. In this case
the necessary conditions $(1)$ to $(3)$ are obvious.

Now we prove the sufficient conditions. Suppose the graphic
sequence $\pi$ satisfies the conditions (1) to (3). Our proof is
by induction on $n$. We first prove the base case where $n=5$. In
this case, $\pi$ is one of the following sequences: $(4^5)$,
$(4^3,3^2)$, $(4^2,3^2,2)$, $(4,3^4)$, $(4,3^3,1)$, $(4,3^2,2^2)$,
$(3^4,2)$, or $(3^3,2,1)$. It is easy to check that all of these
are potentially $(K_5-Y_4)$-graphic. Now suppose the sufficiency
holds for $n-1$ ($n\geq6$), and let $\pi=(d_1,d_2,\cdots,d_n)$ be
a graphic sequence which satisfies (1) to (3). We will prove $\pi$
is potentially
$(K_5-Y_4)$-graphic.\\

\noindent\textbf{Case 1:}  $\pi^\prime=(3^6)$

We have $n=7$ and $\pi$ is one of the following sequences
$(4^3,3^4)$, $(4^2,3^4,2)$ or $(4,3^5,1)$. It is easy to check
that
all of these are potentially $(K_5-Y_4)$-graphic.\\

\noindent\textbf{Case 2:}  $d_n \geq 3$ and $\pi^\prime\neq(3^6)$

Clearly, $d^{\prime}_4 \ge 2$.  If $d_3 \ge 4$, then $d^{\prime}_3
\ge 3$. If $d_3=...=d_n=3$ and $n\ge 6$, $d^{\prime}_3 \ge 3$. It
follows conditions (1) and (2) hold. Thus, by the induction
hypothesis, $\pi^\prime$ is potentially $(K_5-Y_4)$-graphic.
Therefore, $\pi$ is potentially $(K_5-Y_4)$-graphic by Corollary
1.4. In the following, we only consider the cases where
$d_n=2$ or $d_n=1$.\\

\noindent\textbf{Case 3:}  $d_n=2$ and $\pi^\prime\neq(3^6)$

Consider $\pi^\prime=(d_1^\prime,d_2^\prime,\cdots,d_{n-1}^\prime)$.
Since $d_3 \geq 3$ and $d_n = 2$, we have $d_1^{\prime}\geq3$ and
$d_{n-1}^\prime\geq2$. If $d_3^\prime\geq3$, then $\pi^\prime$
satisfies the assumption and it follows $\pi^{\prime}$ is
potentially $(K_5-Y_4)$-graphic. Therefore, $\pi$ is potentially
$(K_5-Y_4)$-graphic by Corollary 1.4. Hence, we may assume
$d_3^{\prime}=2$. We will proceed with the following two cases
$d_1 \geq 4$ and $d_1 = 3$.\\

\noindent\textbf{Subcase 1:} $d_1 \geq 4$

It suffices to consider the case where $d_2=d_3=3$ and
$d_4=d_5=\cdots=d_n=2$. That is, $\pi=(d_1,3^2,2^{n-3})$. Since
$\sigma(\pi)$ is even, $d_1$ must be even. We will prove $\pi$ is
potentially $(K_5-Y_4)$-graphic. It is enough to show
$\pi_1=(d_1-3,2^{n-5},1)$ is graphic. If $d_1 = n-1$, then
$\pi_1=(n-4,2^{n-5},1)$. It suffices to show $\pi_2=(1^{n-5})$ is
graphic. Since $\sigma(\pi_2)$ is even, $\pi_2$ is graphic by
Theorem 1.2. If $d_1\leq n-2$, it suffices to show
$\pi_2=(2^{n-2-d_1},1^{d_1-2})$(or
$\pi_2=(2^{n-1-d_1},1^{d_1-4})$) is graphic. Similarly, one can
show $\pi_2$ is graphic. Thus, $\pi_1=(d_1-3,2^{n-5},1)$ is
graphic
and, in turn, $\pi$ is potentially $(K_5-Y_4)$-graphic.\\

\noindent\textbf{Subcase 2:} $d_1=3$

It suffices to consider the case where $d_1=d_2=d_3=d_4=3$ and
$d_5=\cdots=d_n=2$. That is, $\pi=(3^4,2^{n-4})$.  We will prove
$\pi$ is potentially $(K_5-Y_4)$-graphic. It is enough to show
$\pi_1=(2^{n-5},1^2)$ is graphic. Since $\sigma(\pi_1)$ is even,
$\pi_1$ is graphic by Theorem
1.2 and, in turn, $\pi$ is potentially $(K_5-Y_4)$-graphic.\\

\noindent\textbf{Case 4:} $d_n=1$ and $\pi^\prime\neq(3^6)$

Consider
$\pi^\prime=(d_1^\prime,d_2^\prime,\cdots,d_{n-1}^\prime)$. Since
$d_3 \geq 3$ and $d_4 \geq 2$, we have $d_2^\prime \geq 3$ and
$d_4^{\prime} \geq 2$. If $d_3^{\prime} \geq 3$, then
$\pi^{\prime}$ satisfies the assumptions and it follows
$\pi^{\prime}$ is potentially $(K_5-Y_4)$-graphic. Therefore,
$\pi$ is potentially $(K_5-Y_4)$-graphic by Corollary 1.4. Hence,
we may assume $d_3^\prime=2$. It suffices to consider the case
where $d_1=d_2=d_3=3$ and $d_4=2$. That is,
$\pi=(3^3,2^t,1^{n-3-t})$ where $t \geq 1$ and $n-3-t \geq 1$.
Since $\sigma(\pi)$ is even, $n-t$ must be even. We will prove
$\pi$ is potentially $(K_5-Y_4)$-graphic. It is enough to show
$\pi_1=(2^{t-2},1^{n-2-t})$ is graphic when $t \geq 2$. Since
$\sigma(\pi_1)$ is even, $\pi_1$ is graphic by Theorem 1.2. If
$t=1$, then $\pi=(3^3,2,1^{n-4})$. Similarly we can show
$\pi_2=(1^{n-5})$ is graphic and, in turn, $\pi$ potentially
$(K_5-Y_4)$-graphic.

This completes the proof.\end{proof}

In the remaining of this section, we will use the above two
theorems to find exact values of $\sigma(K_5-P_4,n)$ and
$\sigma(K_5-Y_4,n)$. Note that the value of $\sigma(K_5-P_4,n)$
was determined by Lai in $[9]$ so a much simpler proof is given
here.

\begin{corollary}
(\cite{lai2})
 For $n \geq 5$, $\sigma (K_5-P_4,n)=4n-4$.
\end{corollary}

\begin{proof}
First we claim $\sigma(K_5-P_4,n) \geq 4n-4$ for $n \geq5 $. We
would like to show there exists $\pi_1$ with $\sigma(\pi_1)=4n-6$
such that $\pi_1$ is not potentially $(K_5-P_4)$-graphic. Let
$\pi_1=((n-1)^2,2^{n-2})$. It is easy to see that
$\sigma(\pi_1)=4n-6$ and $\pi_1$ is not potentially
$(K_5-P_4)$-graphic by Theorem 2.1.

Now we show if $\pi$ is an $n$-term $(n\geq5)$ graphical sequence
with $\sigma(\pi) \geq 4n-4$, then there exists a realization of
$\pi$ containing a $K_5-P_4$. If $d_5=1$, then
$\sigma(\pi)=d_1+d_2+d_3+d_4+(n-4)$. Let $X$ be the four vertices
of the largest degrees of $G$ and $Y=V(G)-X$. Since there are at
most six edges in $X$,  $d_1+d_2+d_3+d_4 \leq 12 + |E(X,Y)| \leq
12+(n-4)=n+8$. This leads to $\sigma(\pi) \leq 2n+4 < 4n-4$, a
contradiction. Thus, $d_5 \geq 2$. If $d_2\leq2$, then
$\sigma(\pi) \leq d_1+2(n-1) \leq 3n-3 < 4n-4$, a contradiction.
Thus, $d_2\geq3$. Since $\sigma(\pi) \geq 4n-4$, then $\pi$ is not
one of the following: $(3^2,2^4)$, $(3^2,2^5)$,  and
$(n-1,k,2^t,1^{n-2-t})$ where $n \geq 6$ and $k,t=3,4,\cdots,n-2$,
$(n-k,k+i,2^i,1^{n-i-2})$ where $i=3,4,\cdots,n-2k$ and
$k=1,2,\cdots,[\frac{n-1}{2}]-1$. Thus, $\pi$ satisfies the
conditions (1) to (5) in Theorem 2.1. Therefore, $\pi$ is
potentially $(K_5-P_4)$-graphic by Theorem 2.1. \end{proof}

\begin{corollary}
For $n\geq5$, $\sigma(K_5-Y_4,n)=4n-4$.
\end{corollary}

\begin{proof}
First we claim $\sigma(K_5-Y_4,n) \geq 4n-4$ if $n \geq5$. We
would like to show there exists $\pi_1$ with $\sigma(\pi_1)=4n-6$,
such that $\pi_1$ is not potentially $(K_5-Y_4)$-graphic. Let
$\pi_1=((n-1)^2,2^{n-2})$.  It is easy to see that
$\sigma(\pi_1)=4n-6$ and $\pi_1$ is not potentially
$(K_5-Y_4)$-graphic by Theorem 2.2.

Now we show if $\pi$ is an $n$-term $(n\geq5)$ graphical sequence
with $\sigma(\pi)\geq4n-4$, then there exists a realization of
$\pi$ containing a $K_5-Y_4$.  If $d_4=1$, then
$\sigma(\pi)=d_1+d_2+d_3+(n-3)$. By the similar argument used in
the above corollary, we have $d_1+d_2+d_3 \leq 6+(n-3)=n+3$. This
leads to $\sigma(\pi)\leq2n<4n-4$, a contradiction. Thus,
$d_4\geq2$. Similarly, if $d_3\leq2$, then $\sigma(\pi) \leq
d_1+d_2+2(n-2) \leq2(n-1)+2(n-2)=4n-6<4n-4$, a contradiction.
Thus, $d_3\geq3$. Since $\sigma(\pi)\geq4n-4$, then
$\pi\neq(3^6)$. Thus, $\pi$ satisfies the conditions (1) to (3) in
Theorem 2.2. Therefore, $\pi$ is potentially $(K_5-Y_4)$-graphic
by Theorem 2.2. \end{proof}

\end{document}